\documentclass[12pt]{article}
\usepackage{theorem,amsfonts, amssymb}

\newtheorem{theorem}{Theorem}[section]
\newtheorem{proposition}[theorem]{Proposition}
\newtheorem{lemma}[theorem]{Lemma}
\newtheorem{corollary}[theorem]{Corollary}
\theorembodyfont{\upshape}
\newtheorem{definition}{Definition}[theorem]
\newtheorem{example}[theorem]{Example}
\newtheorem{remark}[theorem]{Remark}
\newtheorem{proof}{Proof}

\newtheorem{acknowledgement}{Acknowledgement}

\newcommand{\bt}{\begin{theorem}}
\newcommand{\et}{\end{theorem}}
\newcommand{\bl}{\begin{lemma}}
\newcommand{\el}{\end{lemma}}
\newcommand{\bp}{\begin{proposition}}
\newcommand{\ep}{\end{proposition}}
\newcommand{\bd}{\begin{definition}}
\newcommand{\bex}{\begin{example}}
\newcommand{\eex}{\end{example}}
\newcommand{\ed}{\end{definition}}
\newcommand{\br}{\begin{remark}}
\newcommand{\er}{\end{remark}}
\newcommand{\bc}{\begin{corollary}}
\newcommand{\ec}{\end{corollary}}
\newcommand{\bo}{\begin{proof}}
\newcommand{\eo}{\end{proof}}
\newcommand{\be}{\begin{enumerate}}
\newcommand{\ee}{\end{enumerate}}

\newcommand{\End}{{\rm End}}

\newcommand{\Aut}{{\rm Aut}}
\newcommand\du{\,{\rm d}}
\newcommand{\Inn}{{\rm Inn}}

\newcommand{\Z}{{\mathbb Z}}
\newcommand{\Q}{{\mathbb Q}}
\newcommand{\N}{{\mathbb N}}
\newcommand{\R}{{\mathbb R}}
\newcommand{\C}{{\mathbb C}}
\newcommand{\T}{{\mathbb T}}
\newcommand{\ceta}{\mbox{\boldmath $\eta$}}

\title{Distal actions and shifted convolution property}
\author{C.\ R.\ E.\ Raja and R.\ Shah }
\date{ }

\begin{document}
\maketitle

\let\epsi=\epsilon
\let\vepsi=\varepsilon
\let\lam=\lambda
\let\Lam=\Lambda 
\let\ap=\alpha
\let\vp=\varphi
\let\ra=\rightarrow
\let\Ra=\Rightarrow 
\let\da=\downarrow 
\let\Llra=\Longleftrightarrow
\let\Lla=\Longleftarrow
\let\lra=\longrightarrow
\let\Lra=\Longrightarrow
\let\ba=\beta
\let\ga=\gamma
\let\Ga=\Gamma
\let\un=\upsilon
\let\mi=\setminus
\let\ol=\overline
\let\ot=\odot

\centerline {To Professor S.\ G.\ Dani on his 60th birthday}

\begin{abstract}
A locally compact group $G$ is said to have shifted convolution 
property (abbr.\ as SCP) if for every regular Borel probability measure 
$\mu$ on $G$, either $\sup _{x\in G} \mu ^n (Cx) \ra 0$ for all compact 
subsets $C$ of $G$, or there exist $x\in G$ and a compact subgroup $K$ 
normalised by $x$ such that $\mu^nx^{-n} \ra \omega_K$, the Haar measure 
on $K$. We first consider distality of factor actions of distal actions. 
It is shown that this holds in particular for factors under compact 
groups invariant under the action and for factors under the connected 
component of identity.  We then characterize groups having 
SCP in terms of a readily verifiable condition on the conjugation action
(point-wise distality).  This has some interesting corollaries to 
distality of certain actions and Choquet Deny measures which actually
motivated SCP and point-wise distal groups.  We also relate distality 
of actions on groups to that of the extensions on the space of probability
measures.
\end{abstract}

\begin{section}{Introduction}

Let $X$ be a Hausdorff space and $\Ga$ be a (topological) semigroup 
acting continuously on $X$ by continuous self-maps.  We say 
that the action of $\Ga$ on $X$ is {\it distal} if for any two distinct 
points $x, y \in X$, the closure  of $\{ (\ap (x), \ap (y)) \mid \ap 
\in \Ga \}$ does not intersect the diagonal $\{ (a, a ) \mid a \in X \}$ 
and we say that the action of $\Ga$ on $X$ is {\it point-wise distal} if 
for each $\gamma\in\Gamma$, $\{\gamma ^n\}_{n\in\N}$-action on $X$ is 
distal.  The notion of distality was introduced by Hilbert 
(cf.\ \cite{El}) and studied by many in different contexts: see 
\cite{El}, \cite{Fu} and \cite{JR} and the references cited therein.

Let $G$ be a locally compact (Hausdorff) group and let $e$ denote the 
identity of $G$.  Let $\Ga$ be a semigroup acting continuously on $G$ 
by endomorphisms.  Then $\Gamma$-action on $G$ is distal if and only if  
$e\not\in\overline {\Gamma x}$ for all $x\in G\mi\{e\}$. The group 
$G$ itself is said to be {\it distal} (resp.\ 
{\it point-wise distal}) if conjugacy action of $G$ on $G$ is distal 
(resp.\ point-wise distal). 

It can easily be seen that the class of distal groups is closed under 
compact extensions.  Abelian groups, discrete groups and 
compact groups are obviously distal. 
Nilpotent groups, connected Lie groups of type $R$ and connected groups of 
polynomial growth are distal (cf.\ \cite{Ro}) and p-adic Lie groups 
of type $R$ and p-adic Lie groups of polynomial growth are point-wise 
distal (cf.\ \cite{R1} and \cite{R2}): point-wise distal groups are called 
non-contracting in \cite{R1} and \cite{Ro}.     

Clearly, distal groups are point-wise distal but 
there are point-wise distal groups which are not distal (see 
\cite {JR} and \cite{Ro} for instance).

For a locally compact group $G$, let $P(G)$ denote the space of all 
regular Borel probability measures on $G$ with weak* topology.  For 
$\mu, \lam \in P(G)$, let $\mu \lam$ denote 
the convolution of $\mu $ and $\lam$ and 
for $x \in G$, let $\delta _x$ denote the Dirac measure at $x$ and $x \mu$ 
(resp.\ $\mu x$) denote the measure $\delta _x \mu$ (resp.\ $\mu 
\delta _x)$.  For $n \geq 1$ and $\mu \in P(G)$, let $\mu ^n$ denote the 
$n$-th convolution power of $\mu$.  For any compact subgroup $K$ of $G$, 
let $\omega_K$ denote the normalized Haar measure on $K$.  For $\mu \in 
P(G)$, let $\check \mu \in P(G)$ be the measure 
defined by $\check \mu (E) = 
\mu (\{ x^{-1} \mid x \in E \})$ for any Borel set $E$ of $G$.   

Let $G$ be a locally compact group and $\Aut(G)$ denote the group of all 
bi-continuous automorphisms of $G$.  If $\Ga $ is a group acting 
continuously on $G$ (by automorphisms), then this action extends to an
action on $P(G)$ which is given by $\ap (\mu ) (E) = \mu (\ap ^{-1} (E))$
for any $\ap \in \Ga$, any $\mu \in P(G)$ and for any measurable subset 
$E$ of $G$.  

We say that a locally compact group $G$ has {\it shifted convolution 
property} which would be called SCP if for any $\mu \in 
P(G)$, one of the following holds:

\be
\item [(i)] the concentration functions of $\mu$, $\sup_{g\in G}\mu ^n(Kg)$
converges to zero (as $n\to\infty$) 
for any compact subset $K$ of $G$ (in which case we say that $\mu$ is 
{\it dissipating}); 

\item [(ii)] there exists a $x \in G$ such that $\mu ^n x^{-n} \ra 
\omega _K$ for some compact subgroup $K$ of $G$ with $x K x^{-1} =K$.  
\ee

If a measure $\mu \in P(G)$ satisfies one of the above two conditions, we 
say that $\mu$ has {\it shifted convolution property} (SCP).  

Tortrat \cite{Tt} proved SCP for groups satisfying certain conditions.  
Motived by \cite{Tt}, Eisele \cite{E2} introduced a class of groups called 
Tortrat groups and 
showed that Tortrat groups have shifted convolution property: a 
locally compact group $G$ is called {\it Tortrat} if for any sequence 
$\{x_n\}$ in $G$ and $\mu \in P(G)$, the sequence $\{x_n \mu x_n^{-1} \}$ 
has idempotent limit point only if $\mu$ is an idempotent.  
Corollary 5.1 and Theorem 5.2 of \cite{E2} showed that the class of 
Tortrat groups strictly contains the class of groups 
satisfying conditions of \cite{Tt}.  Dani and Raja 
\cite{DR2} proved that almost connected groups of polynomial growth, 
equivalently, almost connected (point-wise) distal groups are Tortrat and 
Raja \cite{R4} showed that distal linear groups are Tortrat.  
Thus, proving SCP for almost connected distal groups and distal linear 
groups.  It may be noted that Tortrat groups are distal, but the 
converse is not true in general as we shall see 
in Example \ref{ex}.  Here we show that 
a locally compact group is point-wise distal if and only if it has SCP.  
Our proof/techniques rely on dynamics on compact groups, 
on zero-dimensional groups and on Lie groups over local fields.  It may be 
noted that \cite{B1}-\cite{B3}, \cite{E3}, \cite{J2} and \cite{J1} have 
results on SCP for some classes of groups and measures.  

In Section 2, we state and prove some preliminary results which will
be used often.  In Section 3, we prove some results about factor-actions
of distal actions, where the factors are either modulo compact subgroups
or the connected component of the identity. In Section 4, we prove the main 
result for totally disconnected (metrizable) groups. 
In Section 5, we discuss $\Z$-actions on compact 
metric groups and prove the main 
result for discrete extension of compact groups.  In Section 6, we 
prove the main result for all locally compact groups and 
a few interesting Corollaries.  In Section 7, we compare actions on groups 
and its extension on the space of probability measures. 
\end{section}

\begin{section}{Preliminaries}

The following result is proved in \cite{B1} for locally compact 
$\sigma$-compact groups and is valid without the $\sigma$-compactness 
assumption as we can restrict our attention to the closed subgroup 
generated by the support of $\mu$ which is $\sigma$-compact (see also 
\cite{E1}).  

\bl\label{s}
Let $G$ be a locally compact group and $\mu \in P(G)$.  Suppose 
$\mu$ is non-dissipating.  Then $\mu^n\check\mu^n\ra\rho\in P(G)$ 
and $\mu\rho\check\mu=\rho$. 
\el

We next prove some basic properties of SCP.

\bp\label{bp1}
Let $G$ be a locally compact group and $K$ be a compact 
normal subgroup of $G$.  Let $\pi \colon G\ra G/K$ be the canonical
projection.  If $\mu \in P(G)$ has SCP, then $\pi (\mu )\in P(G/K)$ has 
SCP and hence, if $G$ has SCP, so does $G/K$.   
\ep

\bo
Let $\mu \in P(G)$.  Since $K$ is compact, $\mu$ is dissipating if and 
only if $\pi (\mu )$ is dissipating.  Suppose $\mu$ is non-dissipating 
and 
has SCP.  Then there exists a compact subgroup $L$ of $G$ and $z \in G$ 
such that $\mu ^n z^{-n} \ra \omega _L$ and $z L z^{-1} =L$.   
This implies that 
$\pi (\mu ^n z^{-n}) \ra \omega _{\pi (L)}$ and 
$\pi (z ) \pi (L) \pi (z)^{-1} =\pi (L)$.  
\eo

\bp\label{bp}
Let $G$ be a locally compact group and $\mu \in P(G)$.  Let us assume 
that each neighbourhood $U$ of identity contains a compact normal subgroup 
$K_U$ and $\pi_U \colon G \ra G/K_U$ be the canonical quotient map.  
If $\pi_U(\mu)$ has SCP for all $U$, then $\mu\in P(G)$ also has SCP.  
If $G/K_U$ has SCP for all $U$, then $G$ also has SCP.
\ep

\bo
It is easy to see that for any neighbourhood $U$ of $e$ in $G$, 
$\mu \in P(G)$ is dissipating if and only if $\mu \in P(G/K_U)$ 
is dissipating.  Suppose $\mu \in P(G)$ is 
non-dissipating.  Then $\mu ^n \check \mu ^n \ra \rho \in P(G)$ and 
$\mu \in P(G/K_U)$ is also non-dissipating for all $U$.  If 
$\pi_U(\mu)$ has SCP, then $\{ \pi_U(\mu^n\check\mu ^n)\}$ 
converges to an idempotent normalized by the support of $\pi_U(\mu)$. 
Thus, $\pi_U(\rho)$ is an idempotent normalized by the support of 
$\pi_U(\mu)$.  This 
implies that $\rho \omega_{K_U} = \rho\omega_{K_U}\rho\omega_{K_U}$ 
and $x\rho x^{-1}\omega _{K_U} = \rho\omega _{K_U}$ for all 
$x$ in the support of $\mu$.  Since $\omega_{K_U} \ra \delta _e$ as $U\ra 
e$, we get $\rho = \rho ^2$ and $x\rho x^{-1}= \rho$ for all $x$ in the 
support of $\mu$ (cf.\ \cite{Hey}, Theorem 1.2.2).  Let $x$ be in 
the support of $\mu$.  Then the equation $\mu ^n \rho \check \mu ^n = 
\rho$ implies that $\{ \mu ^n x^{-n} \}$ is relatively compact.  Suppose 
$\nu$ is a limit point of $\{ \mu ^n x^{-n} \}$.  Then 
$\nu \check \nu = \rho$ and $\nu$ is supported on the support of $\rho$.  
By Lemma 2.1 of \cite{E1}, $\nu = \rho z$ for some $z$ in the support of 
$\nu$ and hence $\nu = \rho$.  Thus, $\mu ^n x^{-n} \ra \rho$ for all $x$ 
in the support of $\mu$. 
\eo

\bp\label{bp2}
Let $G$ and $H$ be locally compact groups and let  
$\phi\colon H\ra G$ be a continuous injection.  
If $\mu \in P(H)$ is non-dissipating and 
$\phi (\mu ) \in P(G)$ has SCP, then $\mu$ also has SCP.  In particular, 
$G$ has SCP implies $H$ has SCP. 
\ep

\bo
Let $\mu \in P(H)$ be non-dissipating.  Then $\phi (\mu )\in P(G)$ is 
non-dissipating.  
Suppose $\phi (\mu )$ has SCP.  Then there exists a $x\in G$ such 
that $\phi (\mu ^n) x^{-n} \ra \omega _K$ for some compact group $K$ with 
$xKx^{-1} = K$.  Since $\mu \in P(H)$ is non-dissipating, by Lemma 
\ref{s}, $\mu ^n \check \mu ^n \ra \rho \in P(H)$, 
but $\phi(\mu )^n \phi(\check\mu )^n \ra \omega _K \in P(G)$.  
Thus, $\phi (\rho ) =\omega _K$.  Since $\phi \colon H \ra G$ is 
injective, its extension $\phi \colon P(H) \ra P(G)$ is also 
injective.  Thus, $\rho$ is an idempotent.  We now prove that 
$y\rho y^{-1} = \rho$ for any $y$ in the support of $\mu$.  
Let $y$ be in the support of $\mu$.  Then $ \phi (y )$ is in the support 
of $\phi (\mu )$ and hence $\phi (y) \omega _K \phi (y)^{-1} = \omega _K$.  
Since $\phi $ is injective, $y \rho y^{-1} = \rho$.  We can now show 
as in Proposition 2.2 that $\mu ^n y^{-n} \ra \rho$.
\eo

For a locally compact group $G$ and $\ap\in\Aut(G)$, we define $C_K(\ap)$,
the $K$-{\it contraction group of} $\ap$, for any compact
group $K$ such that $\ap(K)=K$, as follows:
$$
C_K(\ap)=\{x\in G\mid\ap^n(x)K\to K\hbox{ in }G/K\hbox{ as }n\to\infty\}.$$
We denote $C_{\{e\}}(\ap)$ by $C(\ap )$ and it is called the 
{\it contraction group of} $\ap$.  For any
$g\in G$, we denote by $C(g)$, the contraction group of  
$\Inn(g)$, where $\Inn(g)$ is the inner automorphism of $G$
defined by $g$.  An automorphism $\ap$ of a locally 
compact group $G$ is said to contract $G$ if $C(\ap ) = G$.  

The following lemma about contraction groups will be quite useful. 

\bl\label{v1}
Let $G$ be a locally compact group and $\ap \in \Aut(G)$.  Let 
$C(\ap)$ be a non-trivial group.  
Suppose there is a topology on $C(\ap)$ which turns $C(\ap)$ 
into a locally compact 
group $\tilde C(\ap )$ such that $\ap $ is a bi-continuous 
automorphism of $\tilde C(\ap )$ with $\ap ^n (x) \ra e$ as $n \ra \infty$ 
for all $x \in \tilde C(\ap )$ and the canonical inclusion is a continuous 
map of $\tilde C(\ap)$ into $G$.  Then for every probability measure 
$\lam \not = \delta _e$ on $C(\ap )$ such that $\lam$ is supported on a 
compact set of $\tilde C(\ap)$,  $\{\prod_{i=0}^n \ap^i (\lam)\}_{n\in\N}$ 
is convergent in $P(G)$ and its limit is not an idempotent invariant under 
$\ap$. (Here, $\ap^0=I$, the identity map). 
\el

When $G$ is a real or a p-adic Lie group, locally compact topology on 
$C(\ap )$ are possible (cf.\ \cite{Si2}) in which case the Lemma can be 
proved using Theorem 5 of \cite{DS} -- this case is sufficient for our 
purpose.  
Siebert \cite{Si1} has results on topologies of $C(\ap )$ that turn 
$C(\ap )$ into a locally compact group and cases where such 
topologies are not possible.  

\bo $\!\!\!\!\!$ {\bf of Lemma \ref{v1}}\ \  
The convergence of the product $\{\prod_{i=0}^n \ap^i (\lam)\}_{n\in\N}$ 
to, say, $\nu$ in $\tilde C(\ap )$ and hence in $G$ follows from 
``(3) $\Ra$ (2)'' of Proposition 4.3 of \cite{J2}.  This implies that 
$\nu=\lam \ap(\nu )$.  Since $\lam\ne\delta_e$, $\nu\ne\delta_e$. 
Moreover, $\tilde C(\ap)$ does not have any nontrivial compact subgroup 
invariant under $\ap$.  Hence $\nu$ is not an idempotent invariant under
$\alpha$.  
\eo

\br
Given $\ap$ and $\tilde C(\ap)$ as in Lemma \ref{v1} above, there exist
measures $\lam$ on $C(\ap)$ such that 
$\{\lam_n=\prod_{i=0}^n\ap^n(\lam)\}$ is
convergent, but its limit is not an idempotent. For any compactly 
supported nontrivial measure $\lam$ on $\tilde C(\ap)$, 
$\{\lam_n\}$ is convergent and for its limit
point $\nu$, we have $\nu=\lam\ap(\nu)$ and if it were an 
idempotent, say $\omega_K$ for some compact group $K$, then $\ap 
(K)\subset K$ and $\lam$ is supported on $K$. Here, $K\ne\{e\}$ since  
$\lambda\ne\delta_e$. In this case, $\ap (K)$ is strictly contained 
in $K$ and we take $\lam '$ to be any measure on $K$ such that 
$\lam'=\lam'\omega_{\ap (K)}\ne\omega_{\ap(K)}$ and 
$\lam '\not = \omega _K$.  Then $\prod_{i=0}^n\ap^n(\lam') = \lam'$,
for all $n\in\N$. 
\er

The following type of groups occurs in our study often.  Let $\ap$ be an 
automorphism of a locally compact group $G$.  Then 
$\Z \ltimes_\ap G$ denotes the semidirect product of $\Z$ with $G$ where 
the $\Z$-action is given by $\ap$ and is a locally compact group with $G$ 
as an open subgroup. 

\bp\label{lg} 
Let $G$ be a real Lie group.  If $G$ has SCP, then $G$ is point-wise 
distal. In particular, if $\ap \in \Aut(G)$ is such that 
$\Z\ltimes_\alpha G$ has SCP, then $\{ \ap ^n \}_{n \in \Z}$-action 
is distal on $G$.  
\ep

\bo The second assertion trivially follows from the first.
Suppose $G$ is a real Lie group with SCP.  Let $g \in G$.  It is easy to 
see that $C(g)\subset G^0$, the connected component of $e$ in $G$ as 
$G/G^0$ is discrete.  Thus, $C(g)$ has a topology that turns $C(g)$ into 
a locally compact group $\tilde C(g)$ and the identity 
map from $\tilde C(g)\to C(g)$ is continuous. Also, 
conjugacy action of $g$ contracts $\tilde C(g)$. It is obvious that 
$\{g^n \mid n \in \Z \}\cap  C(g) = \{e\}$.  Let $\ap $ denote the inner 
automorphism defined by $g$ restricted to $C(g)$.  This implies that 
there is a continuous injection from $\Z\ltimes _\ap\tilde C(g)$ into 
$G$ given by $(n, h) \mapsto hg^n$.  By, Proposition 
\ref{bp2}, $\Z \ltimes _\ap \tilde C(g)$ has SCP.  

If $C(g)= \tilde C(g)\ne\{e\}$, then by Lemma \ref{v1}, there 
exists a $\lam \in P(\tilde C(g))$ such that 
$e$ is in the support of $\lam$ and $\{\lam_n = \prod _{i=0}^n
\ap ^i(\lam)\}$ converges but its limit is not an idempotent invariant 
under $\ap$.  Now for $\mu=\lam \delta_{1} \in P(\Z \ltimes\tilde C(g))$, 
$\mu^n \delta_{-n}=\lam_{n-1}$ converges but its limit is not an 
idempotent invariant under $\ap$.  This is a contradiction to the above 
assertion that $\Z \ltimes _\ap \tilde C(g)$ has SCP 
(cf.\ \cite{E1}, Theorem 4.3).  Thus, $C(g) = \{e\}$.  
In a similar way one can show that 
$C(g^{-1}) = \{ e \}$.  This implies that the inner automorphism of $g$ is 
distal on the Lie algebra of $G$ (cf.\ \cite{CG}).  Now, Theorem 1.1 of 
\cite{Ab} implies that the conjugacy action of $G$ on $G^0$ is point-wise 
distal.  Since $G/G^0$ is discrete, we get the result.  
\eo
\end{section}

\begin{section}{Factor actions of distal actions}

The following simple result about the factor actions of distal
actions modulo compact groups will be very useful. 

\bt \label{cptf} 
Let $G$ be a Hausdorff group and let $\Gamma$ be a semigroup of 
bi-continuous automorphisms of $G$.  Suppose $K$ is a compact 
subgroup of $G$ such that $\gamma(K)=K$, for all $\gamma\in\Gamma$.  
Then the following are equivalent:
\be

\item  [(1)] $\Gamma$-action on $G$ is distal; 

\item [(2)] $\Gamma$-actions on both $K$ and $G/K$ are distal. 
\ee
\et

Note that if $G$ is compact and metrizable then the above follows from 
Theorem 3.3 of Furstenberg \cite{Fu}. 

\bo
It is sufficient to prove the only non-trivial implication that 
(1) $\Ra$ (2).  Suppose $\Gamma$-action on $G$ is distal.  
Then $\Gamma$-action on $K$ is distal. Let $E$ be the closure of 
$\Gamma$ in $K^K$, the set of
all functions from $K$ into $K$. Then $E$ is a group 
(cf.\ \cite{El}, Theorem 1) and it is a compact subset of $K^K$. 

Let $\{U_d\}$ be a neighbourhood basis at $e$ in $G$ such that 
each $U_d$ is $K$-invariant, i.e.\ $kU_dk^{-1}=U_d$, for all $k\in K$; 
this is possible since $K$ is a compact group.  Suppose $(aK,aK)$ is in 
the closure of $\{(\gamma(x)K,\gamma(y)K)\mid \gamma\in\Gamma\}$ for some 
$x, y \in G$.  Then for any 
neighbourhood $U_d$ of $e$, $\gamma_d(x)K, \gamma_d(y)K\in 
aU_dK=aKU_d$. This implies that 
$\gamma_d(x)=au_dk_d$ and $\gamma_d(y) = au'_dk'_d$
for $u_d, u'_d \in U_d$ and $k_d, k'_d \in K$.  Now 
$\gamma_d(x)k_d^{-1}=au_d$, since $K$ is compact, 
$\{k_d\}$ (resp.\ $\{k'_d\}$) is relatively compact and there 
exists a subnet (we denote it by the same) such that $k_d\to k$ (resp.\ 
$k'_d\to k'$); subnet is possible as $G$ is a topological group.  
Here, $u_d\to e$.  Let $\gamma$ be a limit point of 
$\gamma_d$ in $K^K$. Let 
$k_1=\gamma^{-1}(k^{-1})$ and $k'_1=\gamma^{-1}((k')^{-1})$.
Then passing to a subnet, we get that $\gamma_d(k_1)$ converges to
$\gamma(k_1)=k^{-1}$ and $\gamma_d(xk_1)=au_dk_d\gamma_d(k_1)$,
which converges to $akk^{-1}=a$. Similarly, $\gamma_d(yk'_1)$ converges to 
$a$, therefore $xk_1=yk'_1$ and hence $xK=yK$. 
\eo

The next result follows easily from the above.  

\bc \label{Kcont}
Let $G$ be a locally compact group and let $\ap \in\Aut(G)$ 
be such that $\{ \ap ^n \}_{n \in \N}$-action is distal on $G$. Then for 
any compact subgroup $K$ such that $\ap(K)=K$, we have 
$$C_K^w(\ap) = \{ x\in G \mid {\rm along~a~subsequence}~\ap^n(xK)~{\rm 
converges~to}~K \} = K.$$ 
\ec

\bt \label{connf}
Let $G$ be a locally compact group and let $\Gamma$ be a 
semigroup in $\Aut(G)$. Let $G^0$ denote the connected component of 
the identity $e$ in $G$. Then $\Gamma$-action on $G$ is distal if  and 
only if $\Gamma$-action on both $G^0$ and $G/G^0$ is distal.
\et

\bo
The proof of the ``if'' statement is obvious. Now we prove the ''only if''
statement. Assume that $\Gamma$-action on $G$ is distal. It is obvious 
that the $\Gamma$-action on $G^0$ is distal. We
need to prove that the $\Gamma$-action on $G/G^0$ is distal.

\medskip
\noindent{\bf Step 1:} Since $G^0$ is a connected group it contains a maximal 
compact normal subgroup $C$ such that $G^0/C$ is a real Lie 
group (cf.\ \cite{MZ}) and since $C$ is characteristic, $\gamma(C)=C$ for
all $\gamma\in\Gamma$ and $\Gamma$-action on $G/C$ is 
distal (cf.\ Theorem \ref{cptf}). Also $G/G^0$ is isomorphic to 
$(G/C)/(G^0/C)$.
Hence, without loss of any generality, we may assume that $G^0$ is a 
real Lie group and it has no nontrivial compact normal subgroup. Let $H$ be an 
open subgroup of $G$ containing $G^0$ such that $H/G^0$ is a compact open 
subgroup of $G/G^0$.  Then $H$ contains a compact normal subgroup 
$K$ such that $H/K$ is a real Lie group.  This implies that $KG^0$ is an 
open subgroup of $H$.  Moreover,  
$K\cap G^0$ being a compact normal subgroup of $G^0$, is trivial. 
Thus, $KG^0 \simeq K\times G^0$. In particular, $K\subset Z(G^0)$, the 
centraliser of $G^0$ in $G$, which is a closed normal subgroup. 

\medskip 
\noindent{\bf Step 2}: If $\Gamma$-action on $G/G^0$ is
not distal, then we show that the $\Gamma$-action on $G$ is not distal.
This would lead to a contradiction.  Suppose the identity of $G/G^0$ 
belongs to the closure of $\Gamma(x)G^0$ for some 
$x\not\in G^0$. Let $\{\gamma_d(x)G^0\}$ in $G/G^0$ 
is a net converging to the identity in $G/G^0$. 

We show that there exists an element $k$ in a 
compact totally disconnected subgroup $L$ centralising $G^0$ such that 
$L\cap G^0=\{e\}$, $kG^0=xG^0$ and $\{\gamma_d(k)\}$ converges to $e$.

There exists a neighbourhood basis $\{U_d=K_d\times U'_d\}$, where 
$\{ K_d \}$ is a basis consisting of compact open subgroups contained in 
$K$ and $U'_d$ is a neighbourhood of identity in $G^0$.
We may assume that $\gamma_d(x)\in U_dG^0=K_d\times G^0$. 
Let $\gamma_d(x)=k_dg_d=g_dk_d$ for some $k_d\in K_d$ and $g_d\in G^0$. 

Choose $\gamma$ to be some fixed $\gamma_d$, then $\gamma(x)\in KG^0$. 
Let $\gamma(x)=k'g'$, $k'\in K$ $g'\in G^0$. 
Then $x= kg$, where 
$k=\gamma^{-1}(k')\in\gamma^{-1}(K)\subset Z(G^0)$ 
and $g=\gamma^{-1}(g')\in G^0$ as $G^0$ and $Z(G^0)$ are 
characteristic. Let $L=\gamma^{-1}(K)$ which is a compact
totally disconnected subgroup in $Z(G^0)$ and $L\cap G^0=\{e\}$
and $k\in L$. 

It is enough to show that $\gamma_d(k)\to e$. In fact we show that 
$\gamma_d(k)=k_d\in K_d$ for each $d$.  

$$
\gamma_d(x)= \gamma_d(k)\gamma_d(g)=\gamma_d(g)\gamma_d(k)
=k_dg_d=g_dk_d.$$
Let $a_d= \gamma_d(k)k_d^{-1}=\gamma_d(g^{-1})g_d$.  
Then $a_d\in Z(G^0)\cap G^0$, the centre of $G^0$.
In particular, this implies $k_d$ and 
$a_d$, and hence, $k_d$ and $\gamma_d(k)$ commute.  
Therefore, $a_d$ generates a compact group, which is contained in 
the center of $G^0$, and hence, is trivial, i.e.\ $a_d=e$ and hence,
$\gamma _d(k)=k_d\to e$. This completes the proof.
\eo

The following corollary follows easily from the above theorem. 

\bc
Let $G$ be a locally compact group. Then $G$ is (point-wise) 
distal if and only if $G/G^0$ is (point-wise) distal and the $G$-action on 
$G^0$ is (point-wise) distal.
\ec
\end{section}

\begin{section}{Totally disconnected groups}

In this section we apply Poisson boundary and Choquet-Deny Theorem to show 
that SCP implies point-wise distal for totally disconnected groups.  A 
probability measure $\mu$ on a locally compact group $G$ is said to be 
a Choquet-Deny measure if the bounded continuous functions satisfying the 
equation 
 \begin{equation}
f(g) = \int f(gh) \du\mu (h), ~~~ g\in G 
 \label{eqno1}
 \end{equation}
are constants on 
the cosets of the smallest closed subgroup generated by the support of 
$\mu$.  Let $H_\mu$ denote the space of bounded functions 
satisfying the equation  \begin{equation}
f(g) = \int f(gh) \du\mu (h), ~~~ g\in G
 \label{eqno2}
 \end{equation}
with $L^\infty$-norm.  
If $G$ is a locally compact second countable group, then there exists a 
(compact metric) $G$-space $X$ with a $\sigma$-finite quasi-invariant 
measure 
$\nu$ and an equivariant isometry $\Phi$ of $L^\infty (X, \nu)$ onto 
$H_\mu$ given by the Poisson formula $$\Phi (f) (g) = \int f(gx) \du\rho 
(x) $$ where $\rho$ is a $\mu$-stationary probability measure 
on $X$, that is, $\mu *\rho= \rho$ (cf.\ \cite{J0}).  The $G$-space $X$ is 
called the $\mu$-boundary: see \cite{J0} and \cite{JR} for further details 
on $\mu$-boundary and Choquet-Deny Theorem.  It may be noted that $\mu$ 
is Choquet-Deny if and only if the $\mu$-boundary is trivial.  

Jaworski et al \cite{JRW} proved that any Choquet-Deny measure has SCP.
Here we prove the converse which will be useful in characterizing 
groups with SCP. 

\bp\label{cd} 
Let $G$ be a locally compact group.  If $\mu \in P(G)$ is non-dissipating 
and has shifted convolution property, then $\mu$ is a Choquet-Deny 
measure.
\ep

\bo Let $\mu \in P(G)$ is non-dissipating and has SCP. 
Let $G_\mu$ be the closed subgroup generated by the support of $\mu$.  
We may assume that $G_\mu$ is not compact (cf.\ \cite{JR}, Lemma 4.1).  
Then $\mu \in P(G_\mu )$ is also non-dissipating.  
By Proposition \ref{bp2}, $\mu \in P(G_\mu )$ also has shifted 
convolution property.  
Hence, using Lemma 4.1 of \cite{JR} and replacing $G$ by $G_\mu$, we may 
assume that $G$ is the closed subgroup generated by the support of $\mu$.  

Let us first consider the case when $G$ is second countable.  Since 
$\mu$ has SCP, there is a compact subgroup $K$ and a $z \in G$ such 
that $\mu^n z^{-n}\ra \omega _K$ with $K=zKz^{-1}$.  Then  $\mu$ is 
supported on $zK= Kz$.  This implies that $G/K \simeq {\mathbb Z}$ as $G$ 
is not compact.  Thus, $G$ 
is the semidirect product of $\mathbb Z$ and $K$ where $\mathbb Z$-action 
is given by the inner automorphism of $z$ restricted to $K$.  Also, there 
is a probability measure $\lam$ on $K$ such that $\mu = \lam \delta _z$

By Theorem 4.2 of [JR], the boundary of $\mu$ is a homogeneous space
of $G$ and $K$ acts transitively on the boundary of $\mu$.  Thus, there 
exists a closed subgroup $H$ of $G$ such that $G/H$ is the boundary of 
$\mu$ and $G=KH= HK$ as $K$ is normal in $G$.  Under the natural 
isomorphism between $G/H$ and $K/(K\cap H)$, the canonical action of $G$ on 
$G/H$ defines an action of $G$ on $K/(K\cap H)$ by $g\ot x(K\cap H) = 
ahxh^{-1}(K\cap H)$ where $g = ah$ for $a\in K$ and $h \in H$.  
Let $T\colon K \ra K$ be defined by $T (x) = zxz^{-1}$ for all $x \in K$.  
Let $\eta$ be the Poisson kernel in $K/(K\cap H)$ and 
$\pi\colon K \ra K/(K\cap H)$ be the canonical projection.  
Let $\rho \in P(K)$ be such that $\pi(\rho ) = \eta$ and for 
$n \geq 1$, let $z^n = a_n h_n$ for some $a _n \in K$ and $h _n \in H$.  
Then $\eta = \mu ^n \ot \eta = \mu ^n z^{-n} a_n h_n \ot 
\eta = \pi ( \mu ^n z^{-n} T^n(\rho )a_n) $ for all $n \geq 1$.  Since 
$\mu ^n z^{-n} \ra \omega _K$, it is easy to see that 
$\mu ^n z^{-n} T^n(\rho ) a_n\ra \omega _K$ and hence 
$\eta = \pi (\omega _K)$.  This implies that $\eta$ is $G$-invariant and 
hence the boundary is trivial.  This shows that $\mu$ is a Choquet-Deny 
measure.  

Suppose $G$ is any locally compact group.  Since $G$ is the closed 
subgroup generated by the support of $\mu$, $G$ is $\sigma$-compact.  
Thus, by Theorem 8.7 of \cite{HR}, every neighbourhood $U$ of $e$ contains 
a compact normal subgroup $K_U$ such that $G/K_U$ is second countable.  By 
Proposition \ref{bp1}, $\mu$ is non-dissipating and has shifted 
convolution 
property in $G/K_U$ and hence $\mu$ is a Choquet-Deny measure in $G/K_U$.  
Now, Lemma 3.1 of \cite{JR} implies that $\mu$ is a Choquet-Deny measure 
in $G$.
\eo

\bt\label{tdc}
Let $G$ be a locally compact metrizable group and $G^0$ be the 
connected component of identity in $G$.  Suppose $G$ has SCP.  Then 
$G/G^0$ is point-wise distal.  In particular, a 
totally disconnected locally compact metrizable group with SCP  
is point-wise distal. 
\et

\bo
Since $G^0$ is a connected group it contains a maximal compact normal 
characteristic subgroup $C$ such that $G^0/C$ is a real Lie group 
(cf.\ \cite{MZ}) and $G/C$ also has SCP (cf.\ Proposition \ref{bp1}).  
Also, $G/G^0$ is isomorphic to $(G/C)/(G^0/C)$.  Hence, 
we may assume that $G^0$ is a real Lie group 
and $G^0$ has no nontrivial compact normal subgroup.  
Suppose $G/G^0$ is not point-wise distal. Then there exists 
$g\in G$ such that $\{ g^n \}_{n\in \N}$-action on $G/G^0$ is 
not distal. Let $\ap=\Inn(g)$. Then by 
Proposition 2.1 of \cite{JR}, there exists $x\in G$, such that 
$\ap^n(x)G^0\to G^0$. Now by Step 2 of Theorem \ref{connf}, for 
$\gamma_n=\ap^n$, there exists $k\in L$, a compact totally 
disconnected group in $Z(G^0)$ such that $L\cap G^0=\emptyset$ and 
$\ap^n(k)\to e$. In particular, $C(\ap)=C(g)\ne\{e\}$. 

Now we show that 
$\overline{C(g)}\cap G^0=\{e\}$.  Since $\{g^n\}_{n\in\N}$-action on 
$G/G^0$ is not distal, $\{g^n\}_{n\in\Z}$ generates a discrete infinite 
cyclic group modulo $G^0$, hence, $Z\ltimes_\ap G^0$ is a closed 
subgroup of $G$  and it has SCP 
(cf.\ Proposition \ref{bp2}). Also, since
$G^0$ is a Lie group, by Proposition \ref{lg}, 
$C(g)\cap G^0=\{e\}$. Let $K$ be a subgroup of $G$ as in the 
proof of Theorem \ref{connf}, such that $KG^0=K\times G^0$ and $KG^0$ is 
an open subgroup of $G$.

First we show that $(K\times G^0)\cap C(g)=K\cap C(g)$.  
Let $b\in C(g)$ be such that $b=k'x$, for some $k'\in K$ and 
$x\in G^0$. Then $g^nbg^{-n}\to e$. Arguing as in Step 2 of the proof 
of Theorem \ref{connf}, we get that $g^nk'g^{-n}\to e$ and hence 
$g^nxg^{-n}\to e$. The above implies that 
$x=e$ and $b=k'\in K$. Let $b_n\in C(g)$ be such that $b_n\to b\in G^0$. 
Since $K\times G^0$ is an open neighbourhood of $b$, 
$b_n\in K\times G^0$ for all large $n$.  Since $b_n \in C(g)$, $b_n\in K$ 
and hence $b\in K\cap G^0=\{e\}$.  Thus, $\overline{C(g)}\cap G^0=\{e\}$. 

Now we have that $\overline{C(g)}$ is totally disconnected subgroup 
normalised by $g$. Also, the subgroup 
${\mathbb Z}\ltimes_\ap \overline{C(g)}$ is a (Borel) subgroup of $G$.  
Clearly, the group ${\mathbb Z}\ltimes _\ap \overline{C(g)}$ is totally 
disconnected.  By Lemma 3.5 of \cite{JR}, the measure 
$\mu = \delta_g \nu$ is not a Choquet-Deny measure where 
$\nu  = p\delta_x +(1-p) \delta_e$ for any $0<p<1$ and 
$p\ne {1\over 2}$ for some $x \in C(g)\setminus \{e\}$.  By Lemma 3.3 of 
\cite{JR}, $\mu$ is non-dissipating and hence by 
Proposition \ref{cd}, $\mu$ does not have SCP in ${\mathbb Z}\ltimes_\ap 
\overline{C(g)}$.  By Lemma 3.4 of \cite{JR}, the map $(n, x) \mapsto 
xg^n$ is (continuous) injective from ${\mathbb Z}\ltimes_\ap
\overline{C(g)}$ into $G$ and hence $\mu$ does not have SCP in $G$ also 
(cf.\ Proposition \ref{bp2}).  
\eo

\end{section}

\begin{section} {$\Z$-actions on compact groups}

In this section we apply dynamics of $\Z$-actions on compact groups to 
prove the main theorem for $\Z \ltimes K$.  We first recall some dynamical 
notions.  

Let $K$ be a compact group and $\ap \in \Aut(K)$.  Then  $\ap$ is 
said to be {\it ergodic} (on $K$) if for any $\ap$-invariant Borel 
set $E$ of $K$, $\omega _K(E) = 0 ~~{\rm or }~~1$. We say that 
the $\{ \ap ^n \}_{n\in \Z}$-action on $K$ has 
{\it descending chain condition} which would be called DCC if 
any decreasing sequence $\{K_n\}$ of closed $\ap$-invariant subgroups is 
finite, that is, there exists a $k\geq 1$ such that $K_n = K_m$ for all 
$n, m \geq k$; see \cite{KS} and \cite{S} for details on DCC.  

\bl\label{ks1}
Let $L$ be any non-trivial compact group and $\tau \colon L^\Z \ra L^\Z$ 
be the shift automorphism given by $\tau ((g_i)) = (g_{i+1})$ for any 
$(g_i)\in L^\Z$.  Then the group ${\mathbb Z}\ltimes _\tau L^\Z$ does not 
have SCP.
\el

The above Lemma follows from Theorem \ref{tdc} if $L$ is not connected.  
Here, we give a simple proof in the general case. 

\bo $\!\!\!\!\!$ {\bf of Lemma \ref{ks1}}\ \ 
Let $M = \{(x_i)\in L^\Z \mid x_i = e {\rm ~~ for } ~~ i > 0 \}$.  Then 
$\tau (M) \subset M$.  Let $\mu =\omega _M \delta _1$.  Then 
$\mu ^n \delta _{-n} =  \prod _{i=0} ^{n-1} \tau ^i(\omega _M) 
= \omega _M$ for any $n \geq 1$.  Thus, $\mu$ is non-dissipating and 
$\mu ^n \check \mu ^n \ra\omega _M$.  Suppose $\mu$ has shifted 
convolution property.  Then there exists a compact 
subgroup $N$ of ${\mathbb Z}\ltimes _\tau L^\Z$ and 
$g \in {\mathbb Z}\ltimes _\tau L^\Z$ such that 
$\mu ^n g^{-n} \ra \omega _N$ with $g Ng^{-1} =N$.  Then 
$\mu ^n \check \mu ^n \ra \omega _N$ and hence $M= N$.  Since $g$ 
normalizes $N$, the closed subgroup generated by $g$ and $M$ contains the 
support of $\mu$ and hence $M$ is normalized by the support of $\mu$.  
Since $1$ is in the support of $\mu$, $M$ is $\tau$-invariant.  This is a 
contradiction unless $L$ is trivial.
\eo

\bl\label{ks2}
Let $K$ be a compact metrizable group and $\ap$ be an ergodic 
automorphism of $K$.  Suppose $\Z \ltimes _\ap K$ has SCP and 
the $\{ \ap ^n \}_{n\in \Z}$-action on $K$ has DCC. Then $K$ is 
a compact connected abelian group of finite-dimension.
\el

\bo By Lemma 5.1, $K$ contains no $\ap$-invariant subgroup of the
form $L^{\mathbb Z}$, where $\ap$-action is given by a shift. Now 
using this fact, the assertion essentially follows from the proof 
of Theorem 10.6 of \cite{S}, but for the sake of clarity, 
we give a sketch of proof along the same lines. 

Since $\{ \ap ^n \}_{n\in \Z}$-action on $K$ has DCC, there exists 
a compact Lie group $G$ such that $K\subset G^\Z$ and the action of $\ap$
is given by shift on $G^\Z$ (cf.\ \cite{KS}). For $k\in Z$, let 
$F_H(k)=\{g_k: (g_n)_{n\in\Z}\in K, g_0=e\}$ be closed subgroups of $G$. 
Let $\eta_k:G\to G_k=G/F_H(k)$ be the natural projection and let 
$\ceta_k:K\to (G_k)^\Z$ be the map induced by $\eta_k$ as follows: 
$\ceta_k((g_n)_{n\in\Z})=(\eta_k(g_n))_{n\in\Z}$. Then by the proof of
Theorem 10.6 of \cite{S} (see also \cite{MT} and \cite{Y}), there 
exists $N\geq 1$ such that $\eta_k(G)=G/F_H(k)$ is isomorphic to $\T^n$, 
the $n$-dimensional torus, for all $k\geq N$, for some fixed $n$ and 
$\ceta_N(K)$ is isomorphic to $X^A$, for some $A\in GL_n(\Q)$, where 
$X^A=\{(\pi(y_n) ): y_{n+1}=A(y_n), n\in\Z\}$ and 
$\pi:\R^n\to \R^n/\Z^n\simeq \T^n$ is the canonical quotient map.  Hence 
$\ceta _N(K)$ is a connected abelian finite dimensional group. 

It remains to
show that $\ceta_N$ is an isomorphism, i.e.\ $\ker\ceta_N$ is trivial.
Let $Y_k=\ker\ceta_k$, $k\in\N$. Then $Y_1=\Lambda_1^\Z$ is a subgroup 
of $K$ (cf.\ Proposition 10.2 of \cite{S}), which has SCP.  
By Lemma \ref{ks1} $\Lambda _1$ is trivial and hence, $Y_1$ is trivial.
This shows that $\ceta_1$ is injective, and hence, an isomorphism.  
Now $\ceta_1(Y_2)=\Lambda_2^\Z$ (cf.\ \cite{KS}, Proposition 5.7 (2))
which is a closed subgroup of $\ceta_1(K)$.  Thus, 
$\ceta_1(Y_2)=\Lambda_2^\Z$ has SCP and hence by Lemma \ref{ks1}, 
$\ceta _1(Y_2)$ is trivial.  This implies that $\ceta_2$ is also 
injective.  Arguing inductively, we get that all $\ceta_k$, $k\in\N$, 
are injective.  This completes the proof. 
\eo  

\bl\label{lac3}
Let $K$ be a non-trivial compact metrizable connected finite-dimensional 
abelian group and $\ap$ be an ergodic automorphism of $K$.  Then 
${\mathbb Z}\ltimes _\ap K$ does not have  SCP.  
\el

\bo
Let $m >0$ be such that ${\mathbb Z}^m \subset \hat K \subset 
{\mathbb Q}^m$.  Let $B_m$ be the dual of ${\mathbb Q}^m$.  Then $K$ is a 
quotient of $B_m$.  It can easily be shown that any automorphism of $K$ 
lifts to an automorphism of $B_m$ (see \cite{Ra}).  We denote the lift 
of $\ap$ also by $\ap$.  

For any prime number $p$, let ${\mathbb Q}_p$ be the field of $p$-adic 
numbers.  Then $\ap$ is in 
$GL_m({\mathbb Q}_p)$ for all $p$ and $\ap \in 
GL_m({\mathbb R})$.  We now show that there exists a prime $p$ for which 
$\ap$ has an eigenvalue of p-adic absolute value different from one or 
$\ap \in GL_m({\mathbb R})$ has an eigenvalue of absolute value different 
from one. Suppose for every prime $p$, eigenvalues of $\ap \in GL _m 
({\mathbb Q}_p)$ are of p-adic absolute value one.  
Let $f$ be the characteristic polynomial of $\ap$.  Then the leading 
coefficient of $f$ is one and all coefficients are rational.  Since the 
coefficients of $f$ are elementary symmetric functions of eigenvalues of 
$\ap$, for any prime $p$, $p$-adic absolute value of the coefficients are 
less than or equal one.  This implies that all coefficients of $f$ are 
integers and the leading coefficient is one.  Thus, eigenvalues of 
$\ap \in GL_m({\mathbb R})$ are algebraic integers.  Since $\ap$ is 
ergodic, no eigenvalue of $\ap$ is a root of unity.  Thus, by a 
classical result of Kronecker (see \cite{Es}), we get that 
$\ap \in GL_m({\mathbb R})$ has an 
eigenvalue of absolute value different from one. 
Since $\Z \ltimes _\ap K$ and $\Z \ltimes _{\ap ^{-1}}K$
are isomorphic, we can assume that an eigenvalue of $\ap \in GL_m(\Q _p)$ 
or $\ap \in GL_m({\mathbb R})$ is of absolute value less than one.  

Let $F$ be ${\mathbb Q}_p$ or $\mathbb R$ such that $\ap \in 
GL_m(F)$ has an eigenvalue of absolute value less than one.  
Since ${\mathbb Q}^m \subset F^m$ and ${\mathbb Q}^m$ is dense, there exists 
a continuous injection $\psi \colon F^m \ra B_m$.   It can easily be 
verified that $\psi \ap = \ap \psi$.  
Now, let $V = \{ v\in F^m \mid \ap ^n (v) \ra 0 ~~{\rm as}~~ n\ra \infty 
\}$.  Then $V$ is a non-trivial closed subspace.  
By Lemma \ref{v1} we get a measure 
$\lam \in P(F^m)$ that is supported on $V$ such that the support of 
$\lam$ contains $0$ and $\{\lam\ap(\lam)\cdots\ap^n(\lam)\}_{n\in\N}$ 
converges but the limit point is not an idempotent invariant under $\ap$.  
Let $\pi \colon B_m \ra K$ be the canonical projection.  Let $ \rho \in 
P(F^m )$ be the limit point of $\{\lam\ap(\lam)\cdots\ap^n(\lam)\}_{n\in\N}$.  

Let $(1, \pi \psi (\lam))$ denote the probability measure on 
${\mathbb Z}\ltimes _\ap K$ defined by 
$$(1, \pi \psi (\lam)) (A\times B) = \delta_1 (A) \pi\psi(\lam) (B)$$ 
for any Borel sets $A$ in $\Z$ and $B$ in $K$. We show that
this measure does not have SCP. If possible, suppose 
$(1, \pi\psi (\lam))$ has SCP. Then $\pi\psi(\rho)$ is an 
idempotent invariant under $\ap$ (cf.\ \cite{E1}, Theorem 4.3).  
Let $L$ be the compact 
subgroup of $K$ such that $\pi \psi (\rho ) = \omega _L$.  
Since $\lam$ is supported on $V$ and $V$ is $\ap$-invariant, 
$\rho (V) =1$.  Since $\pi \psi (V)$ is a Borel subgroup of $K$, 
$\pi \psi \rho (\pi \psi (V)) =1$ and hence 
$\omega _L (\pi \psi (V)) =1$.  Since $\omega _L$ is $L$-invariant, 
$\omega _L (x\pi \psi (V)) =1$ for any $x \in L$.  This implies that 
$x\pi \psi (V)\cap \pi \psi (V) \not = \emptyset$ and hence 
$x \in \pi \psi (V)$ for all $x \in L$, that is, 
$L\subset \pi \psi (V)$.  Since $\ap $ contracts $V$, $\ap$ contracts $L$ 
as well.  Since $L$ is a compact group, $L$ is trivial.  Thus, $\psi 
(\lam)$ is supported on the kernel of $\pi$, say, $M$.  

Let $H = \psi ^{-1} (M)\cap V$.  
Then $H$ is a $\ap$-invariant closed subgroup 
of $V$.  Let $\hat M$ be the dual of $M$.  Then $\hat M \simeq 
\Q ^m /\hat K$ has only elements of finite order.  Let $V'$ be the maximal 
vector subspace contained in $H$.  Then $V'$ is $\ap$-invariant.  
Now $\psi$ restricted to $V'$ defines a continuous homomorphism 
$\hat \psi \colon \hat M \ra V'$.  Since $V'$ has no element of finite 
order, $\hat \psi $ and hence $\psi$ is trivial.  This implies that $H$ 
contains no vector subspace.  Thus, $H$ is compact or discrete.  Since 
$\ap$ contracts $H$, $H$ is trivial.  Since $\lam$ is supported on $H$, 
$\lam = \delta _e$.  Then for all $n\in\N$, $\lam\ap(\lam)\cdots\ap^n(\lam)=\delta_e$ an 
idempotent which is $\ap$-invariant, this leads to a contradiction. Therefore,
$(1,\pi\psi(\lam))$, and hence, ${\mathbb Z}\ltimes_\ap K$ does not have SCP.
\eo

\bt\label{cpt}
Let $K$ be a compact group and $\ap$ be an automorphism of $K$.  Suppose 
the group ${\mathbb Z}\ltimes _\ap K$ has SCP.  Then 
$\{ \ap ^n \}_{n \in \Z}$-action is distal on $K$. 
\et

\bo
Let us first consider the case when $K$ is second countable.  Assume that 
the group ${\mathbb Z}\ltimes_\ap K$ has SCP.  By Proposition 2.1 of 
\cite{Ra}, it is enough to show that $\ap$ is not ergodic on $H$ for any 
non-trivial $\Gamma$-invariant closed subgroup $H$ of $K$.
But since ${\mathbb Z}\ltimes _\ap H$ also has SCP, it is enough to show
this for $H=K$. 

Suppose $\ap$ is ergodic on $K$.  
By Theorem 3.16 of \cite{KS}, there exists a decreasing sequence $\{K_i\}$ 
of closed normal $\ap$-invariant subgroups such that 
$\cap K_i = \{ e \}$ and the action of $\ap$ on $K/K_i$ has DCC.  It is 
easy to see that the action of $\ap$ on $K/K_i$ is ergodic.  By 
Proposition \ref{bp1}, ${\mathbb Z}\ltimes _\ap K/K_i$ has SCP. Thus, 
Lemma \ref{ks2} implies that $K/K_i$ is a compact connected abelian 
group of finite-dimension.  Now by Lemma \ref{lac3}, $K=K_i$.  Since 
$\cap K_i = \{ e \}$, $K$ is trivial.  

Now consider the case when $K$ is not necessarily second countable.  
Suppose ${\mathbb Z}\ltimes _\ap K$ has SCP.  Since ${\mathbb Z}\ltimes 
_\ap K$ is $\sigma$-compact, each neighbourhood $U$ of $e$ in $K$ contains 
a compact normal subgroup $K_U$ of ${\mathbb Z}\ltimes _\ap K$ such that 
$({\mathbb Z}\ltimes _\ap K)/K_U$ is second countable (cf.\ \cite{HR}, 
Theorem 8.7).  This implies that 
$K_U$ is a normal subgroup of $K$ and is $\ap$-invariant.  By Proposition 
\ref{bp1}, ${\mathbb Z}\ltimes _\ap (K/K_U)$ has SCP and hence 
$\{\ap^n\}_{n\in\Z}$-action is distal on $K/K_U$.  Since $U$ is an 
arbitrary neighbourhood of $e$ in $K$, 
$\{ \ap ^n \}_{n \in \Z}$-acition is distal on $K$.
\eo

\end{section}

\begin{section}{Locally compact groups}

\bt\label{mt}
A locally compact group $G$ is point-wise distal if 
and only if the group $G$ has SCP. 
\et

\bo We first assume that $G$ is second countable and hence
metrizable. 
Suppose $G$ is point-wise distal.  Let $\mu\in P(G)$ be non-dissipating.  
It is enough to show that $\mu$ has SCP. Without loss of any generality, 
we may assume that $G=G_\mu$, the closed subgroup generated 
by the support of $\mu$.  Let $N_\mu$ be the smallest closed 
normal subgroup of $G$ such that a coset of it contains the support 
of $\mu$.  
Then there exists $z\in G$, a compact subgroup $H$ of $N_\mu$ and an 
increasing sequence $\{k_n\}$ such that $z$ normalises $H$ and there exists 
a Borel subgroup $N_1 \subset \{g \in G \mid z^{-k_n}gz^{k_n}H\to H \}$ 
with $N_\mu = \overline{ N_1 }$ (cf.\ \cite{J1} or \cite{JRW}).   
Let $\ap \colon G 
\ra G$ be $\ap (g) = z^{-1} g z$ for all $g \in G$.  
By Theorem \ref{cptf}, $\{ \ap ^n\}_{n \in \Z}$-action on $G/H$ is 
distal and hence $N_1 \subset H \subset 
N_\mu$.  Thus, $H = N_\mu$.  Now by Lemma 3.8 of \cite{J1}, 
$\mu^n z^{-n}\to \omega_H$, i.e.\ $\mu$ has SCP. 

Conversely, suppose $G$ has SCP.
It is enough to show that $G/G^0$ is point-wise distal and $G$-action on 
$G^0$ is point-wise distal.  It follows from 
Theorem \ref{tdc} that $G/G^0$ is point-wise distal. Now we show that 
$G$-action on $G^0$ is point-wise distal, which in turn would imply that 
$G$ is point-wise distal.  
Let $K$ be the maximal compact 
normal subgroup of $G^0$ such that $G^0/K$ is a real Lie group.  Then $K$ 
is characteristic in $G^0$ and hence $K$ is a normal subgroup 
of $G$.  By Theorem \ref{cpt} $G$-action on $K$ is point-wise distal.
Hence it is enough to show that the $G$-action on $G^0/K$, or 
equivalently, the $G/K$-action on $G^0/K$ is point-wise distal. Also since
$G$ has $SCP$ so does $G/K$ (cf.\ Proposition \ref{bp1}). Now replacing 
$G$ by $G/K$, we may assume that $G^0$ is a Lie group without any 
nontrivial compact normal subgroup. 
Let $x \in G$ and  $\ap _x \colon G^0 \ra G^0 $ be 
defined by $\ap _x (g) = xgx^{-1}$ for all $g \in G^0$. Here, the closed 
subgroup generated by $xG^0$ is either a discrete or a compact subgroup in 
$G/G^0$.  First suppose that it is discrete. Then the map 
$(n, g) \mapsto gx^n$ defines a continuous injection of 
${\mathbb Z}\ltimes _{\ap _x}G^0$ into $G$.  By Proposition \ref{bp2}, 
the group ${\mathbb Z}\ltimes _{\ap _x} G^0$ has SCP.  It follows from 
Proposition \ref{lg} that $G$-action on $G^0$ is point-wise distal. 
 
Now suppose $xG^0$ generates a relatively compact subgroup in $G/G^0$. 
Let $G_x$ be the closed subgroup generated by $G^0$ and $x$. Then it
is an almost connected locally compact group having SCP. Then 
$G_x^0=G_0$ is a Lie group with SCP and by Proposition \ref{lg}, it 
is point-wise distal.  Therefore, $G^0$ is distal. Now $G_x$, being
a compact extension of $G^0$, is also distal.  This implies that 
$\{ \ap _x ^n \}_{n \in \Z}$-action is distal on $G^0$.  Thus, the 
$G$-action on $G^0$ is point-wise distal. 

Now assume that $G$ is any locally compact group. Suppose $G$ is 
point-wise
distal. Let $\mu\in P(G)$ be non-dissipating.  Let $G_\mu$ be the closed
subgroup generated by the support of $\mu$ in $G$. Then $G_\mu$ is
$\sigma$-compact and it is point-wise distal.  This implies that 
every neighbourhood $U$ of the identity $e$ in $G_\mu$ has a compact 
normal subgroup $K_U$ such that $G_\mu /K_U$ is second countable 
(cf.\ \cite{HR}, Theorem 8.7).  
By Theorem \ref{cptf}, $G_\mu/K_U$ is also point-wise distal and hence it 
has SCP.  This implies that the image of $\mu$ on $G_\mu/K_U$ has SCP.  
By Proposition \ref{bp}, $\mu$ itself has SCP . 

Conversely, assume that $G$ has SCP.
Let $x \in G$ and $g \in G$ be such that $e$ is a limit point of 
$\{ x^n g x^{-n} \mid n \in \Z\}$.  Let $H$ be the closed subgroup 
generated by $x$ and $g$.  Then $H$ is $\sigma$-compact and $H$ also has 
SCP.  Then every neighbourhood $U$ of $e$ in $H$ contains a compact normal 
subgroup $K_U$ of $H$ such that $H/K_U$ is second countable 
(cf.\ \cite{HR}, Theorem 8.7).  
By Proposition \ref{bp1}, $H/K_U$ has SCP and hence 
point-wise distal.  This implies $g \in K_U\subset U$.  
Since $U$ is an arbitrary neighbourhood of $e$, $g=e$.  Thus, 
$G$ is point-wise distal. 
\eo

\br 
It is clear from Theorem \ref{cptf} and the above proof that
point-wise distal implies SCP is valid even if $G$ is any metrizable 
group.  
\er

We now prove a few interesting consequences of Theorem \ref{mt}. 

\bc\label{c1}
Let $G$ be a locally compact point-wise distal group.  
If $x \in G$ and $K$ is a compact group such that $xKx^{-1}\subset K$, 
then $xKx^{-1} = K$.  
\ec

This result is known in case $G$ is a real Lie group or  $G$ is a 
Tortrat group (see Lemma 2.2 of \cite{E2}). 

\bo
Let $\lam = \omega _K$.  Then $\lam x\lam x^{-1} = \lam$.  Thus, 
$\prod _{k=0}^n x^{k}\lam x^{-k} = \lam$.  Let $\mu = \lam x $.  Then 
$\mu ^n x^{-n}= \prod _{k=0}^{n-1} x ^k \lam x^{-k}=\lam$ for all 
$n \geq 1$.  Thus, $\mu$ is non-dissipating.  Since 
$G$ is point-wise distal, by Theorem \ref{mt} we get that 
$\{\mu ^n *\check \mu ^n\}$ converges to an idempotent whose support
(group) is normalized by the support of $\mu$.  But 
$\mu^n*\check\mu^n=\lam\check\lam=\omega_K$ and since,
$x$ is in the support of $\mu$, we have $xKx^{-1} = K$. 
\eo  

We next obtain the Krengel-Lin decomposition for measures on 
point-wise distal 
groups.  Let $G$ be a locally compact group with right Haar measure $m$.  
Let $L^2(G)$ be the Hilbert space of all square integrable functions on 
$G$ with respect to $m$ with norm $|| \cdot ||$.  For $\mu \in P(G)$, 
define $P_\mu \colon L^2 (G) \ra L^2(G)$ by 
$$P_\mu (f) (x) = \int f(xy^{-1}) \du\mu (y)$$ for all $x \in G$.  
Then $P_\mu$ is a contraction on $L^2(G)$.  

We recall that a $\mu \in P(G)$ is 
called {\it adapted} if the closed subgroup generated by the 
support of $\mu$ is all of $G$.  

It is shown in \cite{DL} that adapted $\mu \in P(G)$ is dissipating if and 
only if $L^2 (G) = \{f \in L^2 (G) \mid || P_\mu ^n (f) || \ra 0 \} = 
E_0$, (say).  

When $\mu$ is adapted and is non-dissipating, $E_0\not = L^2(G)$ and in 
this situation Krengel-Lin decomposition is about determining the 
orthogonal complement of $E_0$ in $L^2(G)$ and showing that the
orthogonal complement of $E_0$ in $L^2(G)$ is equal to 
$L^2(G, \Sigma _d)$ where $\Sigma _d$ is the deterministic 
$\sigma$-algebra of $P_\mu$ consisting of Borel sets $A$ in $G$ such that 
for each $n \geq 1$ there exists a Borel set $B_n$ with $P_\mu (1_A) = 
1_{B_n}$: see \cite{R3} and references cited therein.  
We now prove the Krengel-Lin decomposition for measures on point-wise 
distal groups which is a continuation of 
\cite{B1}, \cite{KL} and \cite{R3}.  

\bc\label{c2}
Let $G$ be a non-compact locally compact 
point-wise distal group and $\mu$ be an 
adapted probability measure on $G$.  Suppose $\mu$ is non-dissipating.  
Then there exists a compact normal subgroup $K$ such that:

\be
\item $L^2 (G) = E_0 \oplus L^2 (G, \Sigma _d)$ and the deterministic 
$\sigma$-algebra $\Sigma _d$ is generated by $\{x^n K \mid n \in \Z \}$ 
for any $x$ in the support of $\mu$;

\item $(L^2(G, \Sigma _d ), P_\mu )$ is isomorphic to the bilateral shift 
on $l^2 (\Z )$.
\ee
\ec

\bo
The result follows from Theorem \ref{mt} and Proposition 3.1 of \cite{R3}.
\eo

The following is a generalization of implication (a) $\Ra$ (c) of Theorem 
3.6 of \cite{JR}.

\bc\label{cd1}
Let $G$ be a locally compact group.  If every $\mu \in P(G)$ 
is Choquet-Deny, then $G$ is point-wise distal.
\ec

\bo
Suppose $\mu \in P(G)$ is Choquet-Deny.  Then by Theorem 2.25 of 
\cite{JRW} $\mu$ has SCP.  Now Theorem \ref{mt} 
implies that $G$ is point-wise distal. 
\eo

\end{section}

\begin{section}{Actions on spaces of measures}

In this section we consider actions on connected Lie groups and show that 
distality of group actions carries over to the actions on spaces of 
measures under a certain condition.  For a 
connected Lie group $G$ with Lie algebra $\cal G$, $\Aut (G)$ can be 
realised as a subgroup of $GL({\cal G})$, by identifying each 
automorphism with its derivative on $\cal G$.  By an {\it almost algebraic 
subgroup} of $\Aut (G)$ we shall mean a subgroup of $\Aut (G)$ which is a 
subgroup of finite index or equivalently open subgroup in an algebraic 
subgroup of $GL({\cal G})$, under the identification.  Since subgroups of 
$GL({\cal G})$ have algebraic closure, subgroup of an almost algebraic 
group in $\Aut (G)$ has almost algebraic closure, that is, the smallest 
almost algebraic group contining the subgroup.  

\bt\label{ap}
Let $G$ be a connected Lie group having no compact central subgroup of 
positive dimension.  Let $\Gamma\subset\Aut(G)$ be a subgroup.  Then 
the following are equivalent. 

\be
\item [(1)] $\Gamma$-action on $G$ is point-wise distal;

\item [(2)] $\Gamma$-action on $G$ is distal;

\item [(3)] $\Gamma$-action on $P(G)$ is distal;

\item [(4)] $\Gamma$-action on $P(G)$ is point-wise distal.
\ee
\et

\bo
The implication that (1) $\Ra $ (2) follows from Theorem 1.1 of 
\cite{Ab} and \cite{CG} and that (3) $\Ra $ (4) $\Ra$ (1) 
are easy to verify.  We now prove that (2) $\Ra$ (3).  

Assume that $\Gamma$-action on $G$ is distal.  Let 
${\cal G}$ be the Lie algebra of $G$.  Since $\Ga$-action on $G$ is 
distal, by Theorem 1.1 of \cite{Ab} we get that $\Ga$-action on ${\cal G}$ 
is also distal.  By Theorem 1 of \cite{CG}, almost algebraic closure of 
$\Ga$ is a compact extension of a unipotent subgroup of $GL({\cal G})$.  
Since $G$ has no compact central subgroup of positive dimension, Theorem 1 
of \cite{D1} implies that $\Aut(G)$ is an almost algebraic subgroup of 
$GL({\cal G})$.  It follows that almost algebraic closure of $\Ga$ is 
contained in $\Aut(G)$.  Thus, we may assume that $\Ga$ is a unipotent 
subgroup of $\Aut(G)$.  

Since $\Ga$ is a unipotent subgroup of $\Aut (G)\subset GL({\cal G})$, if 
$\Phi$ defines the conjugacy action of the group 
$\tilde G = \Aut (G)\ltimes G$, then $\Phi (\Ga )$ is also a unipotent 
subgroup of $\Aut (\tilde G)$ and hence its action on $\tilde G$ is 
distal.  Thus, replacing $G$ by $\tilde G$ and considering 
the conjugacy action, we may assume that $\ap (x) =x$ for all $x$ in the 
center of $G$ and all $\ap \in \Ga$.  

Suppose $\{ \ap _n \}$ is in $\Ga$ such that $\lim \ap _n (\mu ) = \lim 
\ap _n (\lam ) $ for some $\mu$ and $\lam$ in $P(G)$.  Let 
${\rm Ad} \colon G \ra GL({\cal G})$ be the adjoint representation and 
$\theta \colon \Ga \ra GL( \End ({\cal G}))$ be defined by 
$$\theta (\ap )(v)=(d\ap )v(d\ap )^{-1}$$ for any $\ap \in \Ga$ and any 
$v \in \End ({\cal G})$.  Then $\theta (\ap ) {\rm Ad}(g) = {\rm Ad} (\ap 
(g))$ for any $\ap \in \Ga$ and any $g \in G$.  Thus, 
$\lim \theta (\ap _n ){\rm Ad}(\lam ) = 
\lim \theta (\ap _n){\rm Ad}(\mu )$.  

Now, by passing to subsequence we may assume by Lemma 2.1 of \cite{DR1} 
that $\theta (\ap _n)(w)$ converges for all $w$ in the support of 
${\rm Ad}( \mu )$ and in the support of ${\rm Ad}(\lam )$. 
Let $W= \{ w\in \End ({\cal G}) \mid 
\{ \theta (\ap _n )(w)) \} ~~{\rm converges } \}$.  Then $W$ is a 
subalgebra of $\End ({\cal G})$ containing the supports of 
${\rm Ad}(\lam )$ and ${\rm Ad}(\mu )$.  By Lemma 2.2 of \cite{DR1}, 
there exist sequences $\{\ba _n\}$ and $\{\ba'_n\}$ such that 
$\{\theta (\ba'_n)\}$ converges and $\theta (\ba_n)(v) =v$ for all $v\in W$ 
and $\theta (\ap_n ) = \theta (\ba'_n \ba _n)$ for all $n \geq 1$.  
Again by passing to a subsequence and replacing $\ba_n$ and $\ba_n$ by 
an element from $\ba_n{\rm Ker}(\theta)$ and $\ba_n {\rm Ker}(\theta)$
respectively,  we may assume that $\{\ba'_n\}$ converges 
and $\ap_n = \ba'_n \ba_n$ with $\theta (\ba_n) (w) =w$ for all $w 
\in W$ and all $n \geq 1$.  This implies that $\lim \ba _n (\mu ) = \lim 
\ba _n (\lam )$ and $\ba _n \in \Ga ' = \{ \ap \in \Ga \mid \theta (\ap 
(v)) = v ~~{\rm for ~~ all }~~ v\in W \}$.  Then $\Ga '$ is a
unipotent subgroup of $\Aut(G)$.  

Let $\tilde H$ be the smallest almost algebraic subgroup containing the 
supports of ${\rm Ad }(\mu )$ and ${\rm Ad} (\lam )$.  Now for 
$\ap \in \Ga '$ and $x \in {\rm Ad}^{-1} (\tilde H)=H$, say, $x^{-1} \ap 
(x) $ is in the center of $G$ and $H$ contains the center of $G$.  Thus, 
using Lemma 3.1 of \cite{DR1} we see that 
the conditions of Proposition 3.2 of \cite{DR1} are verified and hence by 
Proposition 3.2 of \cite{DR1}, there exist sequences $\{T _n\}$ and $\{S_n\}$ 
in $\Ga'$ such that $\{T_n\}$ is relatively compact, $S_n(x) =x $ for all $x$ 
in the supports of $\mu$ and $\lam$ and $\ba_n = T_n S_n$ for all $n \geq 
1$.  If $T$ is a limit point of $\{T_n\}$, then $T(\mu ) = T(\lam)$.  Thus, 
$\lam = \mu$.  This proves that the action of $\Ga$ on $P(G)$ is distal.  
\eo

\bex\label{ex}
We will now show that the assumption on the center of $G$ in Theorem 
\ref{ap} can not relaxed.  Let $\T = \{ z\in \C \mid |z| =1 \}$ be the 
circle group and $K = \T^2$ and $\ap \colon K \ra K$ be given by 
$$\ap (w, z) = (wz, z)$$ for all $(w, z) \in K$.  Then $\ap$ is a 
continuous automorphism of $K$ and let $\Ga$ be the group generated by 
$\ap$.  Since the eigenvalues of $\ap$ are one, the action of $\Ga$ on 
$K$ is distal.  Let $L = \{ (1, z )\in K \mid z \in \T \}$ and $\mu 
\in P(K)$ be invariant under $L$.  We will now show that there exists a 
subsequence $\{k_n\}$ such that $\ap ^{k_n} (\mu ) \ra \omega _K$. 

Since $K$ is monothetic, there exists $x, y$ in $\T$ such that the closed 
subgroup generated by $(x, y)$ is $K$.  This implies that the closed 
subgroup generated by $y$ is the circle group $\T$.  Let $\{k_n\}$ be such 
that $y^{k_n} \ra x$.  By passing to a subsequence of $\{k_n\}$, we may 
assume that $\lim \ap ^{k_n} (\mu ) = \lam\in P(K)$ exists.  Now, 
$(x, y) \lam = \lim (y^{k_n}, y) \ap ^{k_n}(\mu ) = \lim \ap ^{k_n} 
((1,y)\mu ) = \lim \ap ^{k_n} (\mu ) = \lam$ as $\mu$ is invariant under 
$L$. Thus, $\lam$ is $(x,y)$-invariant.  Since $K$ is the closed
subgroup generated by $(x, y)$, we get that $\lam$ is $K$-invariant
and hence $\lam$ is the normalized Haar measure on $K$.  This shows that 
if $\mu$ is invariant 
under $L$, then the closure of the orbit $\Ga (\mu )$ contains $\omega 
_K$.  Thus, the action of $\Ga$ on $P(K)$ is not distal.  This also shows 
that the group $\Z \ltimes _\ap K$ is distal but not Tortrat.
\eex
\end{section}

\begin{acknowledgement}
We would like to thank K.\ Schmidt and M.\ Eiensiedler for fruitful 
suggestions in applying Kronecker's Lemma in Lemma \ref{lac3}. The first author would 
like to thank the School of Mathematics, TIFR, Mumbai, India for its hospitality during 
the special trimester ``Lie Groups--Ergodic theory and Probability measures'' when a 
part of work was done.
\end{acknowledgement}

\bigskip\medskip
\advance\baselineskip by 2pt
\begin{tabular}{ll} 
C.\ R.\ E.\ Raja & \hspace*{1cm}  R.\ Shah \\
Stat-Math Unit & \hspace*{1cm}School of Physical Sciences(SPS)\\
Indian Statistical Institute (ISI) & \hspace*{1cm}Jawaharlal Nehru University(JNU)\\
8th Mile Mysore Road & \hspace*{1cm}New Delhi 110 067, India\\
Bangalore 560 059, India & \hspace*{1cm}rshah@mail.jnu.ac.in\\
creraja@isibang.ac.in & \hspace*{1cm}riddhi@math.tifr.res.in
\end{tabular} 

\end{document}